\documentclass[12pt]{article}
\usepackage{latexsym, amsmath, amsfonts, amscd, amssymb, verbatim}
\usepackage{graphicx,graphics}
\usepackage{xcolor}
\usepackage[sort, numbers]{natbib}
\makeatletter   \renewcommand\@biblabel[1]{#1.}     
\makeatother

\def\tto{\;{\lower 1pt \hbox{$\rightarrow$}}\kern -10pt
	\hbox{\raise 2pt \hbox{$\rightarrow$}}\;}

\def\&&{{&\!\!\!\!}}

\newtheorem{Theorem}{Theorem}[section]

\newtheorem{Lemma}{Lemma}[section]

\newtheorem{Definition}{Definition}[section]

\normalsize

\begin{document}

	\title{\bf Local Error Bounds for Affine Variational Inequalities on Hilbert Spaces}
	
	\medskip
	\author{Y. Lim\footnote{Department of Mathematics, Sungkyunkwan University, Suwon 440-746, South Korea; email: ylim@skku.edu.},\ \, H. N. Tuan\footnote{Department of Science Management and International Relations, Hanoi Pedagogical University 2, Vinh Phuc; email: hoangngoctuan@hpu2.edu.vn; be2yes02@yahoo.com.},\ \,
		N. D. Yen\footnote{Institute of Mathematics, Vietnam Academy of Science and Technology, 18 Hoang Quoc Viet, Hanoi 10307, Vietnam;
			email: ndyen@math.ac.vn.}}
	
	\maketitle
	\date{}
	
	\medskip
	\begin{quote}
		\noindent {\bf Abstract.} This paper gives some results related to the research problem about infinite-dimensional affine variational inequalities raised by N.~D.~Yen and X.~Yang [\textit{Affine variational inequalities on normed spaces}, J. Optim. Theory Appl., 178 (2018), 36--55]. Namely, we obtain local error bounds for affine variational inequalities on Hilbert spaces. To do so, we revisit two fundamental properties of polyhedral mappings. Then, we prove a locally upper Lipschitzian property of the inverse of the residual mapping of the infinite-dimensional affine variational inequality under consideration. Finally, we derive the desired local error bounds from that locally upper Lipschitzian property.

		
	\end{quote}

\section{Introduction}\label{Sect_1}
\textit{Error bounds} play an important role in the convergence analysis of numerical algorithms for variational inequalities and  optimization problems. In finite-dimensional settings, there are many important results on error bounds for variational inequalities over polyhedral convex sets, which include linear complementarity problems as a special subclass; see, e.g.,~\cite{Luo_Tseng_1992,Luo_Tseng_1992a,Mangasarian_Shiau_1986,Mathias_Pang_1990,Pang_1987} and the references therein). Applications of the local error bound established in~\cite[Theorem~2.3]{Luo_Tseng_1992} can be found in the same paper and in many other research works (see, e.g.,~\cite{Tuan14, CLY_2022}). Note that the just cited local error bound has a crucial role in proving the $R$-linear convergence of the iterative sequences in question.   
 
\medskip
As far as we know, \textit{infinite-dimensional affine variational inequalities} (infinite-di\-men\-sio\-nal AVIs) have been considered firstly by Yen and Yang~\cite{Yen_Yang_2018}. The authors have shown that quadratic programming problems and linear fractional vector optimization problems on normed spaces can be studied by using affine variational inequalities. They have also obtained two basic facts about infinite-dimensional affine variational inequalities: the Lagrange multiplier rule and the solution set decomposition. In~\cite[Section~5]{Yen_Yang_2018}, the following research problem has been formulated: ``\textit{Solution existence theorems, solution stability, and local error bounds for infinite-dimensional AVIs, similar to those which have been obtained in}~\cite{Luo_Tseng_1992, GowdaPang_1992a,LeeTamYen_book} \textit{for finite-dimensional AVIs, deserve further investigations.}''

\medskip
Although there are many results about error bounds for finite-dimensional AVIs, to the best of our knowledge, there is no literature dealing with error bounds for AVIs in infinite-dimensional settings.
 
\medskip
By using the natural residual function, Luo and Tseng~\cite{Luo_Tseng_1992} obtained a local error bound for finite-dimensional AVIs. Since their proof relies on some properties of finite-dimensional spaces, it cannot be used in an infinite-dimensional setting. However, there is an important remark at~\cite[p.~47]{Luo_Tseng_1992} that one can get the result by applying the locally upper Lipschitzian property of the inverse of the residual function, which can be deduced from a theorem of Robinson~\cite[Proposition~1]{Robinson1981}. Note that Robinson had obtained the theorem on the basis of the Lipschitzian characterization of convex polyhedra of Walkup and Wets~\cite{Walkup_Wets_1969}. We will employ the just mentioned approach in our study. 

\medskip
This paper mainly aims at establishing local error bounds for affine variational inequalities on Hilbert spaces. The aim will be achieved in three steps. First, by a minimax theorem in Simons~\cite{Simons_98}, we propose a refined version of the proof of Theorem~2.207 in Bonnans and Shapiro~\cite{BoSha_2000}, which is an infinite-dimensional analogue of the result of Walkup and Wets~\cite{Walkup_Wets_1969}, to overcome a possible gap in proving that the domain of a generalized convex polyhedral multifunction is a generalized polyhedral convex set (see Section~\ref{Sect_3} for details). Second, we obtain a locally upper Lipschitzian property of the inverse of the residual mapping of the infinite-dimensional AVI under our consideration. Third, local error bounds for AVIs in a Hilbert space setting are derived from the just mentioned locally upper Lipschitzian property.

\medskip
The paper is organized as follows. Section~\ref{Sect_2} presents some preliminaries and formulates the main results. In Section~\ref{Sect_3}, we revisit two fundamental properties of polyhedral mappings which are addressed in~\cite[Theorem 2.207]{BoSha_2000}. Section~\ref{Sect_4} is devoted to the upper Lipschitzian continuity of the inverse of the residual mapping and error bounds for AVIs on Hilbert spaces. 

\section{Preliminaries and the Main Results}\label{Sect_2}

Let $X$ be a Banach space over $\mathbb{R}$ and $X^*$ is the dual space of $X$. We denote  by $\langle x^*,x\rangle$ value of $x^*\in X^*$ at $x\in X$. By $B(x,\varepsilon)$  we denote the open ball in a normed space $X$ with center $x\in X$ and radius $\varepsilon>0$ and $\bar B_X$ denote the
closed unit ball of $X$. Given some vectors $v^1,\dots,v^s$ in $X$,
we denote by ${\rm span}\{v^1,\dots,v^s\}$ the subspace of $X$
generated by $v^1,\dots,v^s$, that is 
$${\rm span}\{v^1,\dots,v^s\}=
\Big\{\sum_{i=1}^s\lambda_iv^i\,:\, \lambda_i\in \mathbb{R}\ \,{\rm for}\
\,i=1,\dots,s\Big\}.$$ 

For any subsets $A$ and $B$ of $X$, the Hausdorff distance between $A$ and $B$ is the quantity
$$h(A,B)=\max\big\{\sup_{x\in A} d(x,B),\,\sup_{y\in B} d(y,A)
\big\}$$ with $d(x,B)$ (resp., $d(y,A)$) denoting the distance from $x$ to $B$ (resp., the distance from $y$ to $A$).

A multifunction $F$ from $X$ to another normed space
$Y$ is called \textit{locally Lipschitzian} at $x_0\in X$ with modulus
$c>0$, if there exists $\varepsilon>0$ such that
$$h\big(F(x),F(x_0)\big)\le c\|x-x_0\|$$
for every $x\in B(x_0,\varepsilon)$ and it is called
\textit{locally upper Lipschitzian} at $x_0\in X$ with modulus
$c>0$, if there exists $\varepsilon>0$ such that
$$F(x)\subset F(x_0)+c\|x-x_0\| \bar B_Y$$
for every $x\in B(x_0,\varepsilon)$. It is not difficult to see that
if $F$ is locally  Lipschitzian at $x_0\in X$ with modulus $c>0$
then $F$ is locally upper Lipschitzian at $x_0\in X$ with modulus
$c>0$. 

\medskip
Following~\cite[p.~133]{BoSha_2000} (see also~\cite{Yen_Yang_2018, Luan_Yen_2020,Luan_Yao_Yen_2018}), we say that $C\subset X$ is a \textit{polyhedral convex set}  if  there exist $x_i^*\in X^*$ and $\alpha_i\in \mathbb R$, $i=1,\dots,m,$ such that
\begin{equation}\label{poly_conv_set}
	C=\{x\in X\; :\; \langle x^*_i,x\rangle\leq \alpha_i,\ i=1,\dots,m\}.
\end{equation}
We say that $K\subset X$ is a \textit{generalized polyhedral convex set}  if it can be represented as the intersection of a polyhedral convex set and a closed affine subspace of $X$, that is, there exist a closed affine subspace $L\subset X$ and $x_i^*\in X^*,$  $\alpha_i\in {\mathbb R}$, $i=1,\dots,m,$ such that
\begin{equation*}\label{gen_poly_conv_set}
	K=\{x\in L\; :\;\langle x^*_i,x\rangle\leq \alpha_i,\ i=1,\dots,m\}.
\end{equation*}

\medskip
Assume that $M:X\to X^*$ is a bounded linear operator, $q\in X^*$ a vector, and $C\subset X$ a polyhedral convex set defined by \eqref{poly_conv_set}. The problem of finding a vector $\bar x\in C$ satisfying \begin{equation}\label{AVI_ineq}
	\langle M\bar x+q,x-\bar x\rangle\geq 0\quad\; \forall x\in C
\end{equation} is called the \textit{affine  variational inequality} (AVI for brevity) defined by the data set $\{M,q,C\}$. In what follows, the solution set of~\eqref{AVI_ineq} will be denoted by $C^*$. 

\medskip
The next two theorems are the main results of this paper on the problem~\eqref{AVI_ineq}. 

\begin{Theorem}\label{R_inverse_is_upper_Lips_2}
	Suppose that $X$ is a Hilbert space whose dual space $X^*$ is identified with $X$, the constraint set $C$ is polyhedral convex, and the operator $M$ has a closed range. Define the residual mapping $R:X\to X$ of~\eqref{AVI_ineq} by setting
	\begin{equation}\label{R(x)} R(x)=x-P_C(x-Mx-q)\end{equation}
	 for every $x\in X$, where $P_C(u)$ denotes the metric projection of $u\in X$ on $C$ {\rm (see~\cite[Sect.~2, Chapter~I]{Kinder_Stam_80})}. Then there exists a
	constant $c>0$ such that the multifunction $R^{-1}:X\rightrightarrows X$ is locally upper Lipschitzian at every point of $X$ with modulus $c$.
\end{Theorem}

\begin{Theorem}\label{Error_bound_for_AVI} Under the assumptions of Theorem~\ref{R_inverse_is_upper_Lips_2}, there exist constants
	$\varepsilon >0$ and $c>0$ such that
	\begin{equation}\label{local_error_bound}
		\begin{array}{rl}
			d(x,C^*)\le c  \|x- \displaystyle P_C\big(x-M x-q)\|
		\end{array}
	\end{equation}
	for every $x\in X$ satisfying the condition $\|x- \displaystyle P_C\big(x-M x-q)\|\le \varepsilon$.
\end{Theorem}

Theorem~\ref{Error_bound_for_AVI} extends the local error bound in~\cite[Theorem
2.3]{Luo_Tseng_1992} to an infinite-dimensional setting. 

\medskip
To prove the above theorems, we will employ the concepts of polyhedral multifunction and generalized polyhedral multifunction, which are formulated as follows.

\begin{Definition}\label{ddef_gpm} {\rm (See~\cite[Definition~2.206]{BoSha_2000}) A multifunction $F : X\rightrightarrows Y$ between normed spaces is said to be \textit{(generalized) polyhedral} if its graph is a union of finitely many (generalized) polyhedral convex 
		sets, called components, in $X\times Y$. In case there is one component, i.e., ${\rm gph}(F)$ is a (generalized)  polyhedral convex set, we say that $F$ is \textit{(generalized) convex 
		polyhedral}.}
\end{Definition}

As noted in~\cite[p.~141]{BoSha_2000}, if $F : X\rightrightarrows Y$ is a generalized convex polyhedral multifunction, then its graph can be written in the form 
\begin{equation}\label{graph_representation}
{\rm gph}F = \left\{ (x, y) \in L\; :\; \langle x^*_i,x\rangle  + \langle x^*_i,x\rangle\leq b_i,\ i = 1, ... ,p\right\},
\end{equation}
where $L$ is a closed affine subspace of $X\times Y$, $x^*_i\in X^*, y^*_i\in Y^*$, and $b_i\in\mathbb R$, 
$i = 1, ... , p$. Moreover, the affine subspace $L$ can be represented in the form 
\begin{equation}\label{affine_set_representation} L := \left\{ (x, y) \in X \times Y : A_1x + A_2y = z\right\},\end{equation}  
where $(A_1, A_2)$ is a continuous linear mapping from $X \times Y$ onto another normed space $Z$. Then one has 
\begin{equation}\label{F(x)_1}
	F(x)=\left\{y \in Y \, : \, \,A_1 x+A_2y=z,\ 
	\langle x_i^*, x \rangle + \langle y_i^*, y \rangle \leq b_i,\ 
	i=1,2,\dots,p\right\}\quad (\forall x\in X).
\end{equation}

In the sequel, we will need the next lemma.

\begin{Lemma}{\rm (A minimax theorem; see \cite[Theorem 3.1]{Simons_98})}\label{Minimax_Theorem}
	Let $A$ be a nonempty convex subset of a vector space, $B$ be a nonempty
	convex subset of a vector space and $B$ also be a compact Hausdorff  topological  space (i.e., $B$ is equipped with a Hausdorff topology and $B$ is compact in that topology).
	Suppose that  $\varphi:\,A\times B\to\mathbb{R}$  is  convex  on  $A$  and is  concave  and  upper  semicontinuous  on  $B$.
	Then $$\inf_{x\in A}\max_{y\in B}\varphi(x,y) = \max_{y\in B}\inf_{x\in A}\varphi(x,y).$$
\end{Lemma}
      
\section{Properties of Generalized Convex Polyhedral Mappings Revisited}\label{Sect_3}

Based on a generalization of Hoffman's error bounds~\cite{Hoffman_1952} to a Banach space setting due to Ioffe~\cite{Ioffe_79}, Bonnans and Shapiro~\cite[Theorem~2.207]{BoSha_2000} have obtained a fundamental theorem on the relationship between generalized convex polyhedral
property and Lipschitzian property of a multifunction and that the domain of a generalized convex polyhedral multifunction is a generalized polyhedral convex set. However, as far as we understand, there might be a gap in proving the latter result. Namely,~\textit{to apply Proposition~2.176 from~\cite{BoSha_2000} about a pair of dual
problems involving an abstract integration to get a minimax equality related to the function $f(x)$ defined at p.~142 in~\cite{BoSha_2000}, the set $E$ in the proof needs to be metric compact while, in general, $E$ is only weakly* compact.} Since Theorem~2.207 from~\cite{BoSha_2000} is the main tool for obtaining Theorems~\ref{R_inverse_is_upper_Lips_2} and~\ref{Error_bound_for_AVI}, the main results of this paper, we are going to propose a refined version of the proof of Bonnans and Shapiro to overcome the just mentioned possible gap.

\medskip
Using the representation~\eqref{F(x)_1}, which can be given for an arbitrary generalized convex polyhedral multifunction between any two normed spaces, we can restate the theorem in question as follows. 

\begin{Theorem}\label{Hausdorff property of multifunction}  {\rm (See~\cite[Theorem 2.207]{BoSha_2000})} Let $X$, $Y$ and $Z$ be Banach spaces, $A_1: X \rightarrow Z$ and $A_2: Y \rightarrow Z$ be linear continuous mappings with $A_2$ having a closed range. Given vectors $z\in Z$, $x_i^* \in X^*$, $y_i^* \in Y^*$, and constants $b_i \in\mathbb{R}$, where $i=1, 2, \dots, p$, one considers the generalized convex polyhedral multifunction $F:\, X\rightrightarrows Y$ represented as
	\begin{equation}\label{F(x)}
	F(x)=\left\{y \in Y \, : \, \,A_1 x+A_2y=z,\ 
	\langle x_i^*, x \rangle + \langle y_i^*, y \rangle \leq b_i,\ 
	i=1,2,\dots,p\right\}
	\end{equation}
	for all $x\in X$. Then, ${\rm dom}F$ is a generalized polyhedral convex set and there
	exists a constant $c>0$, depending on $A_2$ and $y_i^*,
	i=1,2,\dots,p,$ such that
	\begin{equation}\label{Hausdorff_inequality}
	h\big ( F(x_1),F(x_2)\big)\le c ||x_1-x_2||\quad for\ all\ \, x_1,x_2\in {\rm dom}F.
	\end{equation}
\end{Theorem}	

\noindent {\bf Proof.} Repeating the arguments in the proof of \cite[Theorem 2.207]{BoSha_2000}, we find a constant $c>0$ which depends on $A_2$ and $y_i^*,
i=1,2,\dots,p,$ such that~\eqref{Hausdorff_inequality} holds.

By~\eqref{F(x)} we have~\eqref{graph_representation}, where $L$ is given by~\eqref{affine_set_representation}. In particular, ${\rm gph}F$ is a generalized polyhedral convex set in $X\times Y$. Define a linear operator $\Psi_X$ from $X\times Y $ to $X$ by setting $\Psi_X(x,y)=x$ for all $(x,y)\in X\times Y$. From~\eqref{graph_representation} it follows that ${\rm dom}F =\Psi_X({\rm gph}F)$. Thus, by \cite[Proposition~2.10]{Luan_Yao_Yen_2018}, the closure of ${\rm dom}F$ is a generalized polyhedral convex set in $X$. Hence, to justify the assertions of the theorem, it remains to prove that ${\rm dom}F$ is closed. As in the proof of \cite[Theorem 2.207]{BoSha_2000}, setting \begin{equation}\label{g(x)}
g(x)= \inf_{y\in Y}\bigg[\max\Big \{ \|A_1x+A_2y-z\|,\,\max_{1\le i\le p}(\langle x_i^*,x\rangle+\langle y_i^*,y\rangle-b_i)\Big\}\bigg]\end{equation} for all $x\in X$ and using the closedness of the range of $A_2$, as well as Theorem~2.198 from~\cite{BoSha_2000} (see also \cite[Theorem~3.2]{Luan_Yen_2020}), one can show that \begin{equation}\label{domF}{\rm dom}F =\{x\in X\,:\,g(x)\le 0\}.\end{equation} 
To proceed furthermore, put \begin{equation}\label{function_h} h(y,\,\lambda,\gamma)= \langle \lambda,A_1x+A_2y-z\rangle+\sum_{i=1}^p\gamma_i\big(\langle x_i^*,x\rangle+\langle y_i^*,y\rangle-b_i\big)\end{equation}
for all $(y,\lambda,\gamma)\in  Y\times Z^*\times\mathbb{R}^p$
and $E=\big\{(\lambda,\gamma)\in Z^*\times\mathbb{R}_+^p:\, \|\lambda\|+\sum_{i=1}^p\gamma_i\le1\big\}$, where $\gamma=(\gamma_1,\dots,\gamma_p)$. 

\smallskip
{\sc Claim 1}\ \textit{For every $x\in X$, one has}
\begin{equation}\label{minimax_relation}
 g(x)= \displaystyle\inf_{y\in Y}\max_{(\lambda,\gamma)\in E}h(y,\,\lambda,\gamma).
 \end{equation}

In accordance with~\eqref{g(x)}, the claim will be proved if for each $y\in Y$ we can show that
\begin{equation}\label{equality_of_max}
 \max\Big\{ \|A_1x+A_2y-z\|,\,\max_{1\le i\le p}(\langle x_i^*,x\rangle+\langle y_i^*,y\rangle-b_i)\Big\}=\max_{(\lambda,\gamma)\in E}h(y,\,\lambda,\gamma).\end{equation}
By formula~(2.14) from~\cite{BoSha_2000}, where the supremum is attained, we have $$\|A_1x+A_2y-z\|=\max\big\{\langle \lambda,A_1x+A_2y-z\rangle:\,\lambda\in Z^*,\,\|\lambda\|\le1\big\}.$$
Hence there exists $\lambda_0\in Z^*$ satisfying $\|A_1x+A_2y-z\|=\langle \lambda_0,A_1x+A_2y-z\rangle.$
Next, without loss of generality, we can assume that  $$\langle x_1^*,x\rangle+\langle y_1^*,y\rangle-b_1=\max_{1\le i\le p}(\langle x_i^*,x\rangle+\langle y_i^*,y\rangle-b_i).$$

If $\|A_1x+A_2y-z\|\ge\langle x_1^*,x\rangle+\langle y_1^*,y\rangle-b_1$, then
\begin{equation}\label{inequality_of_max_1}
\begin{array}{rl}
&\max\Big \{ \|A_1x+A_2y-z\|,\,\displaystyle\max_{1\le i\le p}(\langle x_i^*,x\rangle+\langle y_i^*,y\rangle-b_i)\Big\}\\
&=\|A_1x+A_2y-z\|=h(y,\lambda_0,0)\le\displaystyle\max_{(\lambda,\gamma)\in E}h(y,\,\lambda,\gamma).
\end{array}
\end{equation}
In addition, for every $(\lambda,\gamma)\in E$, we get
\begin{equation*}
\begin{array}{rl}
h(y,\,\lambda,\gamma)&=\langle \lambda,A_1x+A_2y-z\rangle+\displaystyle\sum_{i=1}^p\gamma_i\big(\langle x_i^*,x\rangle+\langle y_i^*,y\rangle-b_i\big)\\
&\le\| \lambda\|\,\|A_1x+A_2y-z\|+\displaystyle\sum_{i=1}^p\gamma_i\big(\langle x_i^*,x\rangle+\langle y_i^*,y\rangle-b_i\big)\\
&\le(1-\displaystyle\sum_{i=1}^p\gamma_i)\,\|A_1x+A_2y-z\|+\displaystyle\sum_{i=1}^p\gamma_i\big(\langle x_i^*,x\rangle+\langle y_i^*,y\rangle-b_i\big)\\
&= \|A_1x+A_2y-z\|\\
&\hspace{0.5cm}-\displaystyle\sum_{i=1}^p\gamma_i\bigg[\|A_1x+A_2y-z\|-\Big(\langle x_i^*,x\rangle+\langle y_i^*,y\rangle-b_i\Big)\bigg]\\
&\le \|A_1x+A_2y-z\|,
\end{array}
\end{equation*}
where the second inequality follows from the fact that $\|\lambda\|+\sum_{i=1}^p\gamma_i\le1$.
This yields $$\max_{(\lambda,\gamma)\in E}h(y,\,\lambda,\gamma)\le\max\Big\{ \|A_1x+A_2y-z\|,\,\displaystyle\max_{1\le i\le p}(\langle x_i^*,x\rangle+\langle y_i^*,y\rangle-b_i)\Big\}.$$
Combining this with~\eqref{inequality_of_max_1} yields~\eqref{equality_of_max}.

If $\|A_1x+A_2y-z\|<\langle x_1^*,x\rangle+\langle y_1^*,y\rangle-b_1,$ then
\begin{equation}\label{inequality_of_max_2}
\begin{array}{rl}
&\max\Big \{ \|A_1x+A_2y-z\|,\,\displaystyle\max_{1\le i\le p}(\langle x_i^*,x\rangle+\langle y_i^*,y\rangle-b_i)\Big\}\\
&=\langle x_1^*,x\rangle+\langle y_1^*,y\rangle-b_1=h\big(y,0,(1,\ldots,0)\big)\le\displaystyle\max_{(\lambda,\gamma)\in E}h(y,\,\lambda,\gamma).
\end{array}
\end{equation}
Besides, for every $(\lambda,\gamma)\in E$, it holds that
\begin{equation*}
\begin{array}{rl}
h(y,\,\lambda,\gamma)&=\langle \lambda,A_1x+A_2y-z\rangle+\displaystyle\sum_{i=1}^p\gamma_i\big(\langle x_i^*,x\rangle+\langle y_i^*,y\rangle-b_i\big)\\
&\le\| \lambda\|\,\|A_1x+A_2y-z\|+\gamma_1\Big(\langle x_1^*,x\rangle+\langle y_1^*,y\rangle-b_1\Big)\\
&\hspace{0.5cm}+\displaystyle\sum_{i=2}^p\gamma_i\Big(\langle x_i^*,x\rangle+\langle y_i^*,y\rangle-b_i\Big)\\
&\le\| \lambda\|\|A_1x+A_2y-z\|+\big(1-\|\lambda\|-\displaystyle\sum_{i=2}^p\gamma_i\big)\Big(\langle x_1^*,x\rangle+\langle y_1^*,y\rangle-b_1\Big)\\
&\hspace{0.5cm}+\displaystyle\sum_{i=2}^p\gamma_i\Big(\langle x_i^*,x\rangle+\langle y_i^*,y\rangle-b_i\Big)\\
&=\Big(\langle x_1^*,x\rangle+\langle y_1^*,y\rangle-b_1\Big)\\
&\hspace{0.5cm}-\|\lambda\|\bigg[\Big(\langle x_1^*,x\rangle+\langle y_1^*,y\rangle-b_1\Big)-\|A_1x+A_2y-z\|\bigg] \\
&\hspace{0.5cm}-\displaystyle\sum_{i=2}^p\gamma_i\bigg[\Big(\langle x_1^*,x\rangle+\langle y_1^*,y\rangle-b_1\Big)-\Big(\langle x_i^*,x\rangle+\langle y_i^*,y\rangle-b_i\Big)\bigg]\\
&\le\langle x_1^*,x\rangle+\langle y_1^*,y\rangle-b_1,\\
\end{array}
\end{equation*}
where the second inequality is valid because $\|\lambda\|+\sum_{i=1}^p\gamma_i\le1$ and $$\langle x_1^*,x\rangle+\langle y_1^*,y\rangle-b_1>0.$$
Therefore, $$\max_{(\lambda,\gamma)\in E}h(y,\,\lambda,\gamma)\le\max\Big\{ \|A_1x+A_2y-z\|,\,\max_{1\le i\le p}(\langle x_i^*,x\rangle+\langle y_i^*,y\rangle-b_i)\Big\}.$$
Combining this with~\eqref{inequality_of_max_2}, we obtain~\eqref{equality_of_max}. 

We have thus shown that~\eqref{minimax_relation} holds true.

\smallskip
Without loss of generality, we may assume that the norm in the space $\mathbb{R}^p$ is given by the formula $\|\gamma\|=\sum_{i=1}^p|\gamma_i|.$ Then, the norm in the space $Z^*\times\mathbb{R}^p$ is represented as
$$\|(\lambda,\gamma)\|=\|\lambda\|+\sum_{i=1}^p|\gamma_i|,$$ and we have $E=\bar B\cap (Z^*\times\mathbb{R}_+^p) $, where~$\bar B$ denote the closed unit ball of $Z^*\times\mathbb{R}^p$. Clearly, $E$ is convex. Since $\bar B$ and $Z^*\times\mathbb{R}_+^p$ is weakly* closed, $E$ is weakly* closed.  Moreover, by the Banach-Alaoglu theorem, $\bar B$ is weakly* compact. Hence, $E$ is weakly* compact. Since the weak* topology is Hausdorff, $E$ is a compact Hausdorff topological subspace of $Z^*\times\mathbb{R}^p$, which is equipped with the weak* topology.  In addition, for any $y\in  Y$, $h(y,.,.)$ is a  weakly* continuous linear functional on $Z^*\times\mathbb{R}^p$. So, $h(y,.,.)$, is concave and upper continuous on the compact Hausdorff topological space $E$. Furthermore, for each $(\lambda,\gamma)\in Z^*\times\mathbb{R}^p$, by~\eqref{function_h} we can assert that  $h(.,\lambda,\gamma)$ is an affine function; hence it is convex on $Y$. Therefore, by Claim~1 and Lemma~\ref{Minimax_Theorem} we get   
\begin{equation*}
\begin{array}{rl}
 g(x)&= \displaystyle\inf_{y\in Y}\,\max_{(\lambda,\gamma)\in E}h(y,\,\lambda,\gamma)\\
 &=\displaystyle\max_{(\lambda,\gamma)\in E}\,\inf_{y\in Y}h(y,\,\lambda,\gamma)\\
 &=\displaystyle\max_{(\lambda,\gamma)\in E}\,\inf_{y\in Y}\bigg [ \langle \lambda,A_1x+A_2y-z\rangle+\sum_{i=1}^p\gamma_i(\langle x_i^*,x\rangle+\langle y_i^*,y\rangle-b_i)\bigg].
\end{array}
\end{equation*} 
This implies that
\begin{equation}\label{g(x)_new}
 g(x)
  =\displaystyle\max_{(\lambda,\gamma)\in E}\,\inf_{y\in Y}\bigg [ \langle A_2^*\lambda+\sum_{i=1}^p\gamma_iy^*_i,y\rangle+\big(\langle \lambda,A_1x-z\rangle+\sum_{i=1}^p\gamma_i(\langle x_i^*,x\rangle-b_i)\big)\bigg].
 \end{equation} 
For every $(\lambda,\gamma)\in E$, there are two possibilities:

a) $A_2^*\lambda+\displaystyle\sum_{i=1}^p\gamma_i  y^*_i\ne0$. In this case, we notice that  
 	$$\inf_{y\in Y}\bigg [\langle A_2^*\lambda+\displaystyle\sum_{i=1}^p\gamma_iy^*_i,y\rangle+\big(\langle \lambda,A_1x-z\rangle+\displaystyle\sum_{i=1}^p\gamma_i(\langle x_i^*,x\rangle-b_i)\big)\bigg]=-\infty.$$
 	
b) $A_2^*\lambda+\displaystyle\sum_{i=1}^p\gamma_i  y^*_i=0$. In this case, one has 
   \begin{equation*}
   \begin{array}{rl}
   &\inf\limits_{y\in Y}\big [\langle A_2^*\lambda+\displaystyle\sum_{i=1}^p\gamma_iy^*_i,y\rangle+\big(\langle \lambda,A_1x-z\rangle+\displaystyle\sum_{i=1}^p\gamma_i(\langle x_i^*,x\rangle-b_i)\big)\big]\\
   &=\langle \lambda,A_1x-z\rangle+\displaystyle\sum_{i=1}^p\gamma_i(\langle x_i^*,x\rangle-b_i).
     \end{array}
     \end{equation*}
Hence, setting $E'=\big\{(\lambda,\gamma)\in E:\, A_2^*\lambda+\displaystyle\sum_{i=1}^p\gamma_i  y^*_i=0\big\}$ and
$$\bar h\big(x,\lambda,\gamma\big)= \langle \lambda,A_1x-z\rangle+\sum_{i=1}^p\gamma_i(\langle x_i^*,x\rangle-b_i),$$ we deduce from~\eqref{g(x)_new} that $g(x)=\displaystyle\max_{(\lambda,\gamma)\in E'}\bar h\big(x,\lambda,\gamma\big).$ Therefore, thanks to~\eqref{domF}, we have
\begin{equation}\label{domF_new}{\rm dom}F =\{x\in X:\,g(x)\le 0\}=\bigcap_{(\lambda,\gamma)\in E'}\Big\{x\in X:\,\bar h\big(x,\lambda,\gamma\big)\le 0\Big\}.\end{equation}
Since $\bar h(.,\lambda,\gamma)$ is continuous on $X$,
${\rm dom} F$ is closed. (Note that if $E'=\emptyset$, then $g(x)=-\infty$ for every $x\in X$; hence ${\rm dom}F =X$. This comes in full agreement with~\eqref{domF_new}, because the third set there is $X$.) 

The proof is complete.  $\hfill\square$


\section{Proofs of the Main Results}\label{Sect_4}

The forthcoming proofs of our main results will rely on  
the closedness of the domains of generalized convex polyhedral multifunctions of the form~\eqref{F(x)} and the Lipschitzian property~\eqref{Hausdorff_inequality}, which are given by Theorem~\ref{Hausdorff property of multifunction}.

\medskip
{\bf Proof of Theorem~\ref{R_inverse_is_upper_Lips_2}.} Recall that the residual mapping $R:X\to X$ of the affine  variational inequality~\eqref{AVI_ineq} has been defined in~\eqref{R(x)}. To obtain the locally upper Lipschitzian property of the inverse multifunction $R^{-1}:X\rightrightarrows X$, fix a point $y\in {\rm dom} (R^{-1})$.
Taking an arbitrary vector $x\in R^{-1}(y)$, we have $y=R(x)$. Since  $R(x)=x-P_C(x-Mx-q)$, this yields $x-y=P_C(x-M
x-q)$, i.e.,  $x-y$ is the metric projection of $x-Mx-q$ on $C$. Thus, $x-y$ is the unique solution of the strongly convex quadratic programming problem
$$\min\left\{\|(x-Mx-q)-z\|^2: z\in C\right\},$$ which can be rewritten equivalently as
$$\min\left\{\frac{1}{2}\langle z,z\rangle - \langle x-Mx-q, z\rangle: z\in C\right\}.$$
As the derivative of the objective function at $z$ is $z-(x-Mx-q)$ and $C$ is given by~\eqref{poly_conv_set}, by~\cite[Theorem 2.3]{Yen_Yang_2018} we can assert that the last fact is equivalent
to the inclusion $x-y\in C$ together with the existence of Lagrange multipliers
$\lambda_1\ge 0,\cdots,\lambda_m\ge 0$ such that
$$(x-y)-(x-Mx-q)+\sum_{i=1}^m \lambda_i x^*_i =0,$$
i.e., $y-M x-\displaystyle\sum_{i=1}^m\lambda_ix^*_i=q$,
and $\lambda_i(\langle x^*_i, x-y\rangle -\alpha_i)=0$ for $i=1,\cdots m.$ Put $I=\{1,\dots,m\}$,  $$I_0=\{i\in I\, :\, \langle x_i^*,x-y\rangle= \alpha_i\},$$ and $I_1=I\backslash I_0$. Then, setting $\lambda=(\lambda_1,\dots,\lambda_m)$, we see that the triple $(y, x,\lambda)$ satisfies the system
\begin{equation}\label{Euler_eq5}\begin{cases}
y-M x-\displaystyle\sum_{i=1}^m\lambda_ix^*_i=q,\\
\langle x^*_i,x-y\rangle= \alpha_i, \; \lambda_i\geq 0\ \, \forall i\in I_0,\\
\langle x^*_i,x-y\rangle\leq \alpha_i,\ \;
\lambda_i= 0\ \, \forall i\in I_1.
\end{cases}
\end{equation} (A similar observation has been used in the proof of~\cite[Theorem~4.2]{Yen_Yang_2018}.) 
Clearly,~\eqref{Euler_eq5} can be rewritten as
	\begin{equation*}\label{Euler_eq5_2}\begin{cases}
y-M x-\displaystyle\sum_{i=1}^m\lambda_ix^*_i=q,\\
\langle x^*_i,x-y\rangle\le \alpha_i, \;\langle -x^*_i,x-y\rangle\le -\alpha_i, \; -\lambda_i\leq 0\ \, \forall i\in I_0,\\
\langle x^*_i,x-y\rangle\leq \alpha_i,\ \;
\lambda_i\le 0,\ \, -\lambda_i\le 0\ \,\forall i\in I_1.
\end{cases}
\end{equation*}
This is equivalent to
	\begin{equation}\label{Euler_eq5_2}\begin{cases}
y-M x-\displaystyle\sum_{i=1}^m\lambda_ix^*_i=q,\\
\langle -x^*_i,y\rangle+\langle (x^*_i,0_{\mathbb{R}^m}),(x,\lambda)\rangle\le \alpha_i\ \, \forall i\in I_0,\\
\langle x^*_i,y\rangle+\langle (-x^*_i,0_{\mathbb{R}^m}),(x,\lambda)\rangle\le -\alpha_i\ \, \forall i\in I_0,\\
\langle 0_{X},y\rangle+\big\langle (0_{X},0,\cdots,\underbrace{-1}_\text{$i$-th position},\cdots,0),(x,\lambda)\big\rangle\le 0\ \, \forall i\in I_0,\\
\langle -x^*_i,y\rangle+\langle (x^*_i,0_{\mathbb{R}^m}),(x,\lambda)\rangle\le \alpha_i\ \, \forall i\in I_1,\\
\langle 0_{X},y\rangle+\big\langle (0_{X},0,\cdots,\underbrace{1}_\text{$i$-th position},\cdots,0),(x,\lambda)\big\rangle\le 0\ \, \forall i\in I_1,\\
\langle 0_{X},y\rangle+\big\langle (0_{X},0,\cdots,\underbrace{-1}_\text{$i$-th position},\cdots,0),(x,\lambda)\big\rangle\le 0\ \, \forall i\in I_1,
\end{cases}
\end{equation}
where $0_X$ and $0_{\mathbb{R}^m}$ denote the zero vector of $X$ and $\mathbb{R}^m$ respectively.
Observing that the number of inequalities in \eqref{Euler_eq5_2} is $3m$, we can represent the system as
	\begin{equation}\label{Euler_eq5_3}\begin{cases}
y-M x-\displaystyle\sum_{i=1}^m\lambda_ix^*_i=q,\\
	\langle y^*_i,y\rangle+\langle z^*_i,(x,\lambda)\rangle\le \beta_i\ \,  i=1,\cdots,3m,
	\end{cases}
	\end{equation}
where $y^*_i\in X,z^*_i\in X\times\mathbb{R}^m$ and $\beta_i\in\mathbb{R}$ for every $i=1,\cdots,3m$. 

Summing up all the above, we have proved that, for every pair $(y,x)\in X\times X$, the inclusion
$x\in R^{-1}(y)$ holds if and only if there exist
$\lambda=(\lambda_1,\cdots,\lambda_m)$ and a subset $I_0\subset I$ such that the triple $(y, x,\lambda)$ satisfies the system~\eqref{Euler_eq5_3}.

For each subset $I_0\subset I$, denote by $M_{I_0}$ the set of
all $(y,x,\lambda)\in X\times X\times\mathbb R^m$ fulfilling the conditions in~\eqref{Euler_eq5_3} and by $F_{I_0}$
the multifunction from $X$ to $X\times\mathbb{R}^m$
defined by \begin{equation}\label{F_I0} F_{I_0}(y)=\{(x,\lambda)\,:\, (y,x,\lambda)\in
M_{I_0}\}\quad (\forall y\in X).\end{equation} Then one has
\begin{equation}\label{domain_of_inverse_R} {\rm
dom}\,R^{-1}=\displaystyle{\bigcup_{I_0\subset I}} {\rm dom}\,
F_{I_0}
\end{equation}
and
\begin{equation*}\label{image_of_inverse_R} R^{-1}(y)=\pi_X\Big ({\displaystyle\bigcup_{I_0\subset
		I}} F_{I_0}(y)\Big)\ \; \forall y\in  X,
\end{equation*}
where $\pi_X:X\times\mathbb{R}^m\to X$ is the linear operator given by the formula $\pi_X(x,\lambda)=x$ for every $(x,\lambda)\in X\times\mathbb{R}^m$.
Consider the mapping $T:X\times\mathbb{R}^m\to X$ with
	\begin{equation}\label{Operator_T}T(x,\lambda):=-M x-\sum_{i=1}^m\lambda_i x^*_i\end{equation} for every
	$(x,\lambda)\in X\times\mathbb{R}^m$. Clearly, $T$ is a bounded linear operator and
\begin{equation}\label{sum of subspaces} T\big(X\times\mathbb{R}^m\big)=M(X)+ {\rm span}\{x^*_1,\ldots,x^*_m\}.\end{equation}
On one hand,  the linear subspace $M(X)$ of $X$ is closed because $M$ has a closed range by our assumptions. On the other hand, the linear subspace ${\rm
span}\{x^*_1,\ldots,x^*_m\}$ of $X$ is finite dimensional. Therefore, by~\cite[Theorem 1.42]{Rudin_1991}, from~\eqref{sum of subspaces} it follows that the subspace $T\big(X\times\mathbb{R}^m\big)$ is closed in $X$. This means that $T$ has a closed range. Furthermore, by~\eqref{Euler_eq5_3},~\eqref{F_I0}, and~\eqref{Operator_T} we have
\begin{equation*}
F_{I_0}(y)= \bigg \{(x,\lambda) \in X\times\mathbb{R}^m \,:
\, y+T(x,\lambda)=q,\ \langle y^*_i,y\rangle+ \langle
z_i^*,(x,\lambda)\rangle \le \beta_i,\; i=1,\ldots, 3m \bigg \}
\end{equation*}
for every $y\in X$. So, $F_{I_0}$ satisfies all the conditions imposed on the function $F$ in Theorem~\ref{Hausdorff property of multifunction}. Therefore, we can assert that ${\rm
dom}\, F_{I_0}$ is closed and there exists $c_{I_0}>0$ such that
$$
h\big ( F_{I_0}(y_1),F_{I_0}(y_2)\big)\le c_{I_0} ||y_1-y_2||
$$
 for every $y_1,y_2\in {\rm dom}F_{I_0}$. Since $F_{I_0}(y)$ is closed and convex for every $y\in X$ and $X$ is a Hilbert space, the latter implies that 
\begin{equation}\label{Hausdorff_inequality_of_F_I_0}
F_{I_0}(y_1) \subset F_{I_0}(y_2)+ c_{I_0} ||y_1-y_2||\bar B_{ X \times\mathbb{R}^m}, \, \forall y_1,y_2\in {\rm dom}F_{I_0}.
\end{equation}

From~\eqref{domain_of_inverse_R} and the closedness of ${\rm
	dom}\, F_{I_0}$ for each subset $I_0\subset I$, it follows that ${\rm
	dom}  R^{-1}$ is a closed subset of $X$. Setting $\displaystyle c=\max_{I_0\subset I} c_{I_0}$ and fixing a point $\bar y\in X$, we have
\begin{equation}\label{upper_Lipschitz_inclusion}
R^{-1}(y)\subset R^{-1}(\bar y)+c\|y-\bar y\| \bar B_{ X }\end{equation} for
every $y$ in a neighborhood of $\bar y$. Indeed, if $\bar y\notin {\rm dom}
R^{-1}$, then by the closedness of ${\rm dom}  R^{-1}$ we can
find $\delta>0$ such that $B(\bar y,\delta)\cap{\rm
	dom}  R^{-1}=\emptyset$. Hence, $ R^{-1}(y)=\emptyset$ for
every $y\in B(\bar y,\delta)$, and the inclusion
\eqref{upper_Lipschitz_inclusion} is satisfied. Now, suppose that $\bar y\in {\rm dom}  R^{-1}$. By \eqref{domain_of_inverse_R}, we have $\bar y\in \displaystyle{\bigcup_{I_0\subset
		I}} {\rm dom}\, F_{I_0}$. Denote by $\mathcal{A}(\bar y)$ the family of all subsets $I_0\subset I$ such that
$\bar y\in {\rm dom} F_{I_0}$ and $\mathcal{B}(\bar y)$ the complement of  $\mathcal{A}(\bar y)$ in the family of all subsets of $I$. Since $\displaystyle\bigcup_{I_0\in \mathcal{B}(\bar y)}\, {\rm dom}  F_{I_0}$ is closed, there
exists $\delta>0$ satisfying $ B(\bar y,\delta)\cap
\Big(\displaystyle\bigcup_{I_0\in \mathcal{B}(\bar y)}\, {\rm dom}  F_{I_0}\Big)=\emptyset$. Let $y\in B(\bar y,\delta)$ be given arbitrarily. If
$y\in{\rm dom} F_{I_0}$ for some $I_0\subset I$, then $I_0\in \mathcal{A}(\bar y)$; so $\bar y\in{\rm dom} F_{I_0}$. Therefore,  it follows from~\eqref{Hausdorff_inequality_of_F_I_0} that
$$F_{I_0}(y)\subset F_{I_0}(\bar y)+c\|y-\bar y\|\bar B_{ X \times\mathbb{R}^m}.$$
If $y\notin{\rm dom} F_{I_0}$ then $F_{I_0}(y)=\emptyset$. Hence the latter inclusion is also valid. Thus, 
$$F_{I_0}(y)\subset F_{I_0}(\bar y)+c\|y-\bar y\|\bar B_{ X \times\mathbb{R}^m}$$
for every $I_0\subset I$. Consequently, 
$${\displaystyle\bigcup_{I_0\subset
		I}}F_{I_0}(y)\subset \displaystyle\bigcup_{I_0\subset
	I}F_{I_0}(\bar y)+c\|y-\bar y\|\bar B_{X\times\mathbb{R}^m}.$$
Acting the linear operator $\pi_X$ on both sides of the last inclusion, we get~\eqref{upper_Lipschitz_inclusion}. 

We have thus proved that the multifunction $R^{-1}:X\rightrightarrows X$ is locally upper Lipschitzian at every point of $X$ with modulus $c$. $\hfill\square$

\medskip
We are now in a position to prove Theorem~\ref{Error_bound_for_AVI}, which asserts that the local error bound~\eqref{local_error_bound} holds for some constants $\varepsilon >0$ and $c>0$, provided that the assumptions of Theorem~\ref{R_inverse_is_upper_Lips_2} are fulfilled.

\medskip
{\bf Proof of Theorem~\ref{Error_bound_for_AVI}.} By Theorem~\ref{R_inverse_is_upper_Lips_2}, there exists $c>0$ such that the inverse of the map  $$R(x)=x-P_C(x-Mx-q)$$
is upper Lipschitzian at every point $x\in X$ with the same modulus $c$. In particular, there exists $\varepsilon>0$ such that \begin{equation}\label{inclusion_new}
R^{-1}(z)\subset R^{-1}(0)+c \|z\| \bar B_{X}\end{equation} 
for every
$z\in X$ satisfying $\|z\|\le\varepsilon$. 

We have  $R^{-1}(0)=C^*$. Indeed, a vector $\bar x\in X$ belongs to $R^{-1}(0)$ if and only if $R(\bar x)=0$ or, equivalently, $\bar x=P_C(\bar x-M\bar x-q)$. By the characterization of the metric projection on a closed convex set (see~\cite[Theorem 5.2]{Brezis_2011}), the last equality is equivalent to
$\langle M\bar x+q,x- \bar x\rangle\ge 0$ for $x\in C$.
This means that $\bar x\in C^*$.

Now, fix any vector $x\in X$ with $\|x- \displaystyle
P_C\big(x-Mx-q)\|\le \varepsilon$ and put $z=R(x)$. We have
$$\|z\|=\|R(x)\|=\|x- \displaystyle
P_C\big(x-Mx-q)\|\le\varepsilon.$$ Hence, by~\eqref{inclusion_new} one gets $$ R^{-1}(z)\subset R^{-1}(0)+c \|z\|
\bar B_{X} =C^*+c \|z\|
\bar B_{X}.$$ This implies that $ x\in C^*+c \|z\|
\bar B_{X}$; so we can find $\bar x\in C^*$ such that $ x\in
\bar x+c \|z\|\bar B_{X}$. It follows that $\|x-\bar x\|\le c \|z\|$. Consequently, we obtain
$$d(x,C^*)\le\|x-\bar x\|\le c \|z\|=c \|x- \displaystyle
P_C\big(x-Mx-q)\|.$$ 
 Thus, we have shown that the error bound~\eqref{local_error_bound} holds for every $x\in X$ satisfying the condition $$\|x- \displaystyle P_C\big(x-M x-q)\|\le \varepsilon.$$

The proof is complete. $\hfill\square$


\begin{thebibliography}{10}
	\addcontentsline{toc}{chapter}{\numberline{References}}
	
	\bibitem{BoSha_2000}   Bonnans, J.F., Shapiro, A.:  Perturbation Analysis of Optimization Problems, Springer- Verlag, New York (2000)
	
	\bibitem{Brezis_2011}    Brezis, H.:  Functional Analysis, Sobolev Spaces and Partial Differential Equations, Springer-Velag, New York (2011)
	
	\bibitem{GowdaPang_1992a}   Gowda, M.S., Pang, J.S.: On the boundedness and stability of solutions to the affine variational inequality problem, SIAM J. Control Optim. {\bf 32}, 421--441  (1994) 
	
	\bibitem{CLY_2022} Cuong,  T.H., Lim, Y., Yen, N.D.: On a solution method in indefinite quadratic programming under linear constraints, Optimization, Published online: 03 Nov 2022, https://doi.org/10.1080/02331934.2022.2141056. 
	
	\bibitem{Hoffman_1952}   Hoffman,  A.J.: On approximate solutions of systems of linear inequalities, J. Research Nat. Bur. Standards {\bf 49}, 263--265 (1952)
	
	\bibitem{Ioffe_79}    Ioffe,  A.D.:  Regular points of Lipschitz functions, Trans. Amer.
	Math. Soc. {\bf 251}, 61--69 (1979)
	
	\bibitem{Kinder_Stam_80}    Kinderlehrer, D., Stampacchia, G.: An Introduction to Variational Inequalities and Their Applications, Academic Press, New York (1980)
	
	\bibitem{LeeTamYen_book}     Lee, G.M.,  Tam, N.N., Yen,  N.D.:  Quadratic Programming and Affine Variational Inequalities: A Qualitative Study, Springer Verlag (2005)
	
	\bibitem{Luan_Yao_Yen_2018}    Luan,  N.N., Yao, J.C., Yen,   N.D.: On some generalized polyhedral convex constructions, Numer. Funct. Anal. Optim. {\bf 39}, 537--570 (2018)
	
	\bibitem{Luan_Yen_2020}   Luan, N.N., Yen, N.D.: A representation of generalized convex polyhedra and applications, Optimization {\bf 69}, 471--492 (2020)  
	
	\bibitem{Luo_Tseng_1992}    Luo, Z.Q., Tseng, P.:  Error bound and convergence analysis of matrix splitting algorithms for the affine variational inequality problem, SIAM J. Optim.  {\bf 2}, 43--54 (1992)
	\bibitem{Luo_Tseng_1992a}    Luo, Z.Q., Tseng, P.: On a global error bound for a class of monotone affine variational inequality problems, Oper. Res. Lett. {\bf 11}, 159--165 (1992)
	
	\bibitem{Mangasarian_Shiau_1986}   Mangasarian, O.L., Shiau, T.H.: Error bounds for monotone linear complementarity problems, Math. Programming {\bf 36}, 81--89 (1986)
	
	\bibitem{Mathias_Pang_1990}   Mathias, R., Pang,  J.S.: Error bounds for the linear complementarity problem with a P-matrix, Linear Algebra Appl. {\bf 132}, 123--136 (1990)
	
	\bibitem{Pang_1987} Pang,  J.S.: A posteriori error bounds for the linearly-constrained variational inequality problem, Math. Oper. Res. {\bf 12}, 474--484 (1987)
	
	\bibitem{Robinson1981}    Robinson, S.M.:  Some continuity properties of polyhedral multifunctions, Math. Programming Stud.  {\bf 14}, 206--214 (1981)
	
	\bibitem{Rudin_1991}   Rudin, W.: Functional Analysis, 2nd Edition, McGraw Hill, New York (1991)
	
	\bibitem{Simons_98}  Simons, S.: Minimax and  Monotonicity,  Lecture Notes in  Mathematics  1693, Springer Verlag (1998)
	
	\bibitem{Tuan14}   Tuan,  H.N.: Linear convergence of a type of iterative sequences in nonconvex quadratic programming,  J. Math. Anal. Appl.  {\bf 423}, 1311--1319 (2015)
	
	\bibitem{Walkup_Wets_1969}   Walkup, D.W., Wets,  R.J.B.: A Lipschitzian characterization of convex polyhedra, Proc. Amer. Math. Soc. {\bf 23}, 167--173 (1969) 
	
	\bibitem{Yen_Yang_2018}      Yen,  N.D., Yang,  X.: Affine variational inequalities on normed spaces, J. Optim. Theory Appl.  {\bf 178}, 36--55 (2018)
	
\end{thebibliography}
\end{document}